\documentclass[]{interact}

\usepackage[utf8]{inputenc}

\usepackage{epstopdf}
\usepackage[caption=false]{subfig}

\usepackage[numbers,sort&compress]{natbib}
\bibpunct[, ]{[}{]}{,}{n}{,}{,}
\renewcommand\bibfont{\fontsize{10}{12}\selectfont}

\usepackage[unicode]{hyperref}

\usepackage{latexsym}
\usepackage{indentfirst}
\usepackage{amsxtra}
\usepackage{amssymb}
\usepackage{amsthm}
\usepackage{amsmath}

\usepackage{mathrsfs} 

\usepackage[capitalise]{cleveref}

\newtheorem{theor}{Theorem}
\newtheorem*{theor*}{Theorem}
\newtheorem{prop}[theor]{Proposition}

\newtheorem*{cor*}{Corollary}
\theoremstyle{definition}               

\newtheorem{ex}[theor]{Example}

\DeclareMathOperator{\Sym}{Sym}

\DeclareMathOperator{\id}{id}

\newcommand{\lambdaa}[2]{\lambda_{#1}{#2}}
\newcommand{\lambdaam}[2]{\lambda^{-1}_{#1}{#2}}
\newcommand{\rhoo}[2]{\rho_{#1}{#2}}
\newcommand{\alphaa}[3]{\alpha^{#1}_{#2}{#3}}
\newcommand{\betaa}[3]{\beta^{#1}_{#2}{#3}}

\begin{document}

\articletype{}

\title{The matched product of set-theoretical solutions associated with shelves}

\author{
\name{
Francesco CATINO\textsuperscript{a}\thanks{CONTACT  F. Catino.  Email: francesco.catino@unisalento.it,
Address: Dipartimento di Matematica e Fisica ``Ennio De Giorgi",
Università del Salento -
Via Provinciale Lecce - Arnesano , 73100 Lecce (Italy)}, 
Ilaria COLAZZO\textsuperscript{a}, 
and 
Paola STEFANELLI\textsuperscript{a}
}
\affil{\textsuperscript{a}Dipartimento di Matematica e Fisica ``Ennio De Giorgi",
	Università del Salento, Lecce, Italy}
}

\maketitle

\begin{abstract}
We investigate the matched product of solutions associated with right and left shelves. 
First, we prove that the requirements to provide the matched product of solutions that come from shelves can be simplified. Then we give conditions for left non-degeneracy of the matched product. 
Later, we compute the structure shelf of  the matched product of solutions. Finally, we prove that the structure shelf of the matched product does not depend on the choice of the actions. 
\end{abstract}

\begin{amscode}
	16T25, 81R50, 16Y99, 6N20  
\end{amscode}

\section{Introduction}

Shelves, also known as self-distributive systems, became a subject of systematic study in the 1980s when they gave fundamental computable invariants of braids and knots \cite{De20book, ElNe15book}. Furthermore, there is a strict connection between braid groups $B_n$ and self distributive systems that appeared in the 1980s and led to the discovery in 1991 of a left invariant linear order on the braid groups, see \cite{Deh94} and \cite[Chapter $7$]{KaTu08}.

On the other hand, the study of the Yang-Baxter equation, named after the authors of the first papers in which the equation arose, Yang \cite{Ya67} and Baxter \cite{Ba72},  has been an extensive research area for the past sixty years since is a fundamental tool in several different fields such as statistical mechanics, quantum group theory, and low-dimensional topology.
In 1992, V. Drinfel$'$d \cite{Dr90} suggested focusing on a specific class of solutions: the \emph{set-theoretical solutions} or \emph{braided sets}. Namely, given a set $X$, a set-theoretical solution, shortly a solution, is a map $r: X\times X \to X \times X$ such that the following equation 
\begin{align*}
	\left(r\times \id_X\right)\left(\id_X\times r\right)\left(r\times \id_X\right) = \left(\id_X\times r\right)\left(r\times \id_X\right)\left(\id_X\times r\right)
\end{align*}
is satisfied. Let $r$ be such a solution on X. For $x, y \in X$, define the maps $\lambda_x : X \to X$ and $\rho_y : X \to X$ by $r\left(x, y\right) = \left(\lambdaa{x}{(y)}, \rhoo{y}{(x)}\right)$, for all $x,y \in X$. A solution $r$ is said to be \emph{left} (resp. \emph{right}) \emph{non-degenerate} if $\lambdaa{x}{}$ (resp. $\rhoo{x}{}$) is bijective, for each $x \in X$. 
The seminal papers by Etingof, Schedler, and Soloviev \cite{ESS99}, and Gateva-Ivanova and Van den Bergh \cite{GaB98}, laid the foundations for studying the class of non-degenerate solutions that are also involutive. A solution $r$ on $X$ is said to be involutive if $r^2 = \id_{X\times X}$.
To describe the class of all non-degenerate involutive solutions, Rump introduced two algebraic structures named brace \cite{Ru07a, CeJO14} and cycle set. Many authors attacked the problem of finding constructions of braces, see \cite{Ru07b, CCoSt16, Ga18, Sm18, BaCeJ18, Ba18, CeGaSm18, CeSmVe19} to name a few. Soloviev in \cite{So00} and Lu, Yan, and Zhu in \cite{LuYZ00} studied non-degenerate bijective not necessarily involutive solutions. Guarnieri and Vendramin \cite{GVe17} generalized the concept of braces to obtain an algebraic structure called skew braces in correspondence with bijective solutions. Such solutions have been relatively investigated \cite{SmVe18, Ch18, CCoSt19, ElNeTs19, NZ19}. Further advancements in the field of skew braces relating to Hopf-Galois structures can be found in \cite{SmVe18, Ch18, NZ19}, whereas \cite{CCoSt19}, an extension of \cite{CCoSt15}, partially answered the extension problem in a simplified case. Lebed in \cite{Le17} drew attention on idempotent solutions that, although of little interest in physics, provide a tool for dealing with very different algebraic structures ranging from free (commutative) monoids to factorizable monoids, and from distributive lattices to Young tableaux and plactic monoids. 
Catino, Colazzo, and Stefanelli \cite{CCoSt17}, and Jespers and Van Antwerpen \cite{VAJe19} introduced the algebraic structure called semi-brace to deal with solutions that are not necessarily non-degenerate or that are idempotent. Given a set $X$, a solution $r$ is said to be idempotent if $r^2=r$. In \cite{CCoSt19x} the authors proved that the matched product of solutions, a novel construction technique for solutions of the Yang-Baxter equation introduced in \cite{CCoSt18x}, is a unifying tool for treating solutions with finite order, that includes involutive and idempotent solutions as particular cases.

In the last years, the connection between shelves and the solutions of the Yang-Baxter equation was highlighted. In \cite{LeVe17}, Lebed and Vendramin studied the (co)homology of left non-degenerate solutions. Additionally, they showed how to associate to any left non-degenerate solution a shelf operation that captures many of its properties: invertibility, involutivity, the structure group. In \cite{Le18} Lebed deepened the link between a left non-degenerate solution and its associated shelf.

In this work, we investigate the matched product of solutions associated with right and left shelves. First, we focus on the matched product of left shelves, and we prove that we can simplify the requirements to provide the matched product of solutions related to left shelves. Moreover, we prove that we obtain the left non-degeneracy if and only if we start with left racks. Later, we treat solutions associated with right shelves. Here too, we show that the requirements to obtain a matched product system of solutions are straightforward.
Furthermore, we present the matched product of one solution associated with a left shelf and the other associated with a right shelf; and we obtain natural conditions to check the properties of the matched product system. In all three cases, we prove that the structure shelf of the matched product solution does not depend on the maps actions that give the matched product. Finally, we prove that this is a general result, i.e., that the structure shelf of the matched product of left non-degenerate solutions does not depend on the actions.

\section{Definitions and preliminary results}

At first, we recall the strict connection between shelves (or self-distributive systems) and the solutions of the Yang-Baxter equation. For further details we refer to the extensive treatment by Lebed in \cite{Le18}. 

A set $X$ with an operation $\triangleright$ is a \emph{left shelf} if $\triangleright$ is a left self-distributive operation, i.e., 
$x\triangleright\left(y\triangleright z\right)
= \left(x\triangleright y\right)\triangleright\left(x\triangleright z\right)$, for all $x,y,z\in X$. Likewise, a set $X$ with a right  self-distributive operation $\triangleleft$ (i.e., such that $\left(x\triangleleft y\right)\triangleleft z = \left(x\triangleleft z\right)\triangleleft \left(y\triangleleft z\right)$, for all $x,y,z\in X$) is said to be a \emph{right shelf}.\\
A \emph{left rack} is a left shelf $\left(X, \triangleright\right)$ such that the maps $L_{x}:X\to X$ defined by $L_{x}\left(y\right):= x\triangleright y$, for all $x,y\in X$, are bijections. 
In the same way, a \emph{right rack} is a right shelf $\left(X, \triangleleft\right)$ such that the maps $R_{x}:X\to X$ defined by $R_{y}\left(x\right):= x\triangleleft y$, for all $x,y\in X$, are bijections.

From any shelf it is possible to obtain a set-theoretical solution of the Yang-Baxter equation. Indeed
$\left(X, \triangleright\right)$ is a left shelf if and only if the map $r_{\triangleright}:X\times X\to X\times X$ defined by $r_{\triangleright}\left(x,y\right) = \left(y, y\triangleright x\right)$, for all $x,y\in X$, is a solution. Clearly $r_{\triangleright}$ is left non-degenerate. 
In addition, $\left(X, \triangleright\right)$ is a left rack if and only if the map $r_{\triangleright}$ is a non-degenerate solution.
Analogously, $\left(X, \triangleleft\right)$ is a right shelf if and only if the map $r_{\triangleleft}:X\times X\to X\times X$ defined by $r_{\triangleleft}\left(x,y\right) = \left(y\triangleleft x, x\right)$, for all $x,y\in X$, is a solution and $\left(X, \triangleleft\right)$ is a right rack if and only if the map $r_{\triangleleft}$ is a non-degenerate solution.

On the other hand, if $r$ is a left non-degenerate solution on $X$, $r\left(a,b\right) = \left(\lambda_{a}\left(b\right), \rho_{b}\left(a\right)\right)$, with $\lambda_{a},\rho_{b}$ maps from $X$ into itself, for all $a, b\in X$,
then the binary operation $\triangleright_{r}$ defined by 
\begin{align*}
	a\triangleright_{r} b = \lambdaa{a}{\rhoo{\lambdaam{b}{\left(a\right)}}{\left(b\right)}}
\end{align*} 
gives to $X$ a structure of a shelf called the \emph{structure shelf of $r$}. 
In \cite{LeVe17} it is shown by a diagrammatic proof that $\triangleright_{r}$ is a self-distributive operation while here we prove this result in purely algebraic terms. 
Indeed, it is easy to verify that $r$ is a solution if and only if 
\begin{align}
	&\lambdaa{a}{\lambdaa{b}{}} = \lambdaa{\lambdaa{a}{\left(b\right)}}{\lambdaa{\rhoo{b}{\left(a\right)}}{}}\label{eq:lambdalambda}\\
	&\rhoo{\lambdaa{\rhoo{b}{\left(a\right)}}{\left(c\right)}}{\lambdaa{a}{\left(b\right)}} =\label{eq:lambdarho} \lambdaa{\rhoo{\lambdaa{b}{\left(c\right)}}{\left(a\right)}}{\rhoo{c}{\left(b\right)}}\\
	&\rhoo{c}{\rhoo{b}{}} = \rhoo{\rhoo{c}{\left(b\right)}}{\rhoo{\lambdaa{b}{\left(c\right)}}{}},\label{eq:rhorho}
\end{align}
for all $a,b,c \in X$. 
Since $r$ is left non-degenerate, using \eqref{eq:lambdalambda} we can rewrite \eqref{eq:lambdarho} as:
\begin{align}
	\rhoo{\lambdaam{\lambdaa{a}{\left(b\right)}}{\left(c\right)}}{\lambdaa{a}{\left(b\right)}}
	=\lambdaa{\rhoo{\lambdaam{a}{\left(c\right)}}{\left(a\right)}}{\rhoo{\lambdaam{b}{\lambdaam{a}{\left(c\right)}}}{\left(b\right)}}\label{eq:lambdarho1}\\
	\lambdaam{\rhoo{\lambdaam{a}{\left(c\right)}}{\left(a\right)}}{\rhoo{\lambdaam{b}{\left(c\right)}}{\left(b\right)}}
	=\rhoo{\lambdaam{\lambdaam{a}{\left(b\right)}}{\lambdaam{a}{\left(c\right)}}}{\lambdaam{a}{\left(b\right)}}.\label{eq:lambdarho2}
\end{align}
Hence,
\begin{align*}
	a \triangleright_{r}\left(b\triangleright_{r} c\right)
	&=\lambdaa{a}{\rhoo{\lambdaam{\lambdaa{b}{\rhoo{\lambdaam{c}{\left(b\right)}}{\left(c\right)}}}{\left(a\right)}}{\lambdaa{b}{\rhoo{\lambdaam{c}{\left(b\right)}}{\left(c\right)}}}}\\
	&=\lambdaa{a}{\lambdaa{\rhoo{\lambdaam{b}{\left(a\right)}}{\left(b\right)}}{\rhoo{\lambdaam{\rhoo{\lambdaam{c}{\left(b\right)}}{\left(c\right)}}{\lambdaam{b}{\left(a\right)}}}{\rhoo{\lambdaam{c}{\left(b\right)}}{\left(c\right)}}}}&\mbox{by \eqref{eq:lambdarho1}}\\
	&=\lambdaa{a}{\lambdaa{\rhoo{\lambdaam{b}{\left(a\right)}}{\left(b\right)}}{\rhoo{\lambdaam{\lambdaam{c}{\left(b\right)}}{\lambdaam{c}{\left(a\right)}}}{\rhoo{\lambdaam{c}{\left(b\right)}}{\left(c\right)}}}}&\mbox{by \eqref{eq:lambdalambda}}
\end{align*}
and 
\begin{align*}
	&\left(a\triangleright_{r}b\right)\triangleright_{r}\left(a\triangleright_{r} c\right)=\lambdaa{\lambdaa{a}{\rhoo{\lambdaam{b}{\left(a\right)}}{\left(b\right)}}}{\rhoo{\lambdaam{\lambdaa{a}{\rhoo{\lambdaam{c}{\left(a\right)}}{\left(c\right)}}}{\lambdaa{a}{\rhoo{\lambdaam{b}{\left(a\right)}}{\left(b\right)}}}}{\lambdaa{a}{\rhoo{\lambdaam{c}{\left(a\right)}}{\left(c\right)}}}}\\
	&=\lambdaa{\lambdaa{a}{\rhoo{\lambdaam{b}{\left(a\right)}}{\left(b\right)}}}{\lambdaa{\rhoo{\lambdaam{a}{\lambdaa{a}{\rhoo{\lambdaam{b}{\left(a\right)}}{\left(b\right)}}}}{\left(a\right)}}{\rhoo{\lambdaam{\rhoo{\lambdaam{c}{\left(a\right)}}{\left(c\right)}}{\lambdaam{a}{\lambdaa{a}{\rhoo{\lambdaam{b}{\left(a\right)}}{\left(b\right)}}}}}{\rhoo{\lambdaam{c}{\left(a\right)}}{\left(c\right)}}}}&\mbox{by \eqref{eq:lambdarho1}}\\
	&=\lambdaa{\lambdaa{a}{\rhoo{\lambdaam{b}{\left(a\right)}}{\left(b\right)}}}{\lambdaa{\rhoo{\rhoo{\lambdaam{b}{\left(a\right)}}{\left(b\right)}}{\left(a\right)}}{\rhoo{\lambdaam{\rhoo{\lambdaam{c}{\left(a\right)}}{\left(c\right)}}{\rhoo{\lambdaam{b}{\left(a\right)}}{\left(b\right)}}}{\rhoo{\lambdaam{c}{\left(a\right)}}{\left(c\right)}}}}\\
	&=\lambdaa{a}{\lambdaa{\rhoo{\lambdaam{b}{\left(a\right)}}{\left(b\right)}}{\rhoo{\lambdaam{\rhoo{\lambdaam{c}{\left(a\right)}}{\left(c\right)}}{\rhoo{\lambdaam{b}{\left(a\right)}}{\left(b\right)}}}{\rhoo{\lambdaam{c}{\left(a\right)}}{\left(c\right)}}}}&\mbox{by \eqref{eq:lambdalambda}}\\
	&=\lambdaa{a}{\lambdaa{\rhoo{\lambdaam{b}{\left(a\right)}}{\left(b\right)}}{\rhoo{\rhoo{\lambdaam{\lambdaam{c}{\left(b\right)}}{\lambdaam{c}{\left(a\right)}}}{\lambdaam{c}{\left(b\right)}}}{\rhoo{\lambdaam{c}{\left(a\right)}}{\left(c\right)}}}}&\mbox{by \eqref{eq:lambdarho2}}\\
	&=\lambdaa{a}{\lambdaa{\rhoo{\lambdaam{b}{\left(a\right)}}{\left(b\right)}}{\rhoo{\lambdaam{\lambdaam{c}{\left(b\right)}}{\lambdaam{c}{\left(a\right)}}}{\rhoo{\lambdaam{c}{\left(b\right)}}{\left(c\right)}}}}, &\mbox{by \eqref{eq:rhorho}}
\end{align*}
i.e., $\triangleright_{r}$ is self-distributive. 
Clearly, if $r_{\triangleright}$ is a solution associated with a  left shelf $\left(X,\triangleright\right)$ then the structure shelf $\left(X,\triangleright_{r_{\triangleright}}\right)$ is exactly $\left(X,\triangleright\right)$.

The aim of this work is to analyze the matched product of solutions that come from shelves. For convenience of the reader, here we briefly recall the definition and some basic results of the matched product contained in \cite{CCoSt18x}.
The matched product is a novel construction technique which allows for obtaining new solutions of the Yang-Baxter equation on the cartesian product of sets, starting from absolutely arbitrary solutions.\\
Namely, given a solution $r_S$ on a set $S$ and a solution $r_T$ on a set $T$, if $\alpha: T \to \Sym\left(S\right)$ and $\beta: S \to \Sym\left(T\right)$ are maps, set $\alpha_u:=\alpha\left(u\right)$, for every $u \in T$, and $\beta_a:=\beta\left(a\right)$, for every $a\in S$, then the quadruple
$\left(r_S,r_T, \alpha,\beta\right)$ is said to be a \emph{matched product system of solutions} if the following conditions hold

{\small
	\begin{center}
		\begin{minipage}[b]{.5\textwidth}
			\vspace{-\baselineskip}
			\begin{align}\label{eq:primo}\tag{s1}
				\alpha_u\alpha_v = \alpha_{\lambda_u\left(v\right)}\alpha_{\rho_{v}\left(u\right)}
			\end{align}
		\end{minipage}%
		\hfill\hfill\hfill
		\begin{minipage}[b]{.5\textwidth}
			\vspace{-\baselineskip}
			\begin{align}\label{eq:secondo}\tag{s2}
				\beta_a\beta_b=\beta_{\lambda_a\left(b\right)}\beta_{\rho_b\left(a\right)}
			\end{align}
		\end{minipage}
	\end{center}
	\begin{center}
		\begin{minipage}[b]{.5\textwidth}
			\vspace{-\baselineskip}
			\begin{align}\label{eq:quinto}\tag{s3}
				\rho_{\alpha^{-1}_u\!\left(b\right)}\alpha^{-1}_{\beta_a\left(u\right)}\left(a\right) = \alpha^{-1}_{\beta_{\rho_b\left(a\right)}\beta^{-1}_b\left(u\right)}\rho_b\left(a\right)
			\end{align}
		\end{minipage}%
		\hfill\hfill
		\begin{minipage}[b]{.5\textwidth}
			\vspace{-\baselineskip}
			\begin{align}\label{eq:sesto}\tag{s4}
				\rho_{\beta^{-1}_a\!\left(v\right)}\beta^{-1}_{\alpha_u\left(a\right)}\left(u\right) = \beta^{-1}_{\alpha_{\rho_v\left(u\right)}\alpha^{-1}_v\left(a\right)}\rho_v\left(u\right)
			\end{align}
		\end{minipage}
	\end{center}
	\begin{center}
		\begin{minipage}[b]{.5\textwidth}
			\vspace{-\baselineskip}
			\begin{align}\label{eq:terzo}\tag{s5}
				\lambda_a\alphaa{}{\betaa{-1}{a}{\left(u\right)}}{}= \alphaa{}{u}{\lambdaa{\alphaa{-1}{u}{\left(a\right)}}{}}
			\end{align}
		\end{minipage}%
		\hfill\hfill
		\begin{minipage}[b]{.5\textwidth}
			\vspace{-\baselineskip}
			\begin{align}\label{eq:quarto}\tag{s6}
				\lambdaa{u}{\betaa{}{\alphaa{-1}{u}{\left(a\right)}}{}}=\betaa{}{a}{\lambdaa{\betaa{-1}{a}{\left(u\right)}}{}}
			\end{align}
		\end{minipage}
	\end{center}
}
\noindent for all $a,b \in S$ and $u,v \in T$.

As shown in {\cite[Theorem 1]{CCoSt18x}}, any matched product system of solutions determines a new solution on the set $S\times T$. Namely, if $\left(r_S,r_T, \alpha,\beta\right)$ is a matched product system of solutions, then the map $r:S{ \times} T\times S{ \times} T \to S{ \times} T\times S{ \times} T$ defined by
\begin{align*}
	&r\left(\left(a, u\right), 
	\left(b, v\right)\right) := 
	\left(\left(\alphaa{}{u}{\lambdaa{\bar{a}}{\left(b\right)}},\, \betaa{}{a}{\lambdaa{\bar{u}}{\left(v\right)}}\right),\ \left(\alphaa{-1}{\overline{U}}{\rhoo{\alphaa{}{\bar{u}}{\left(b\right)}}{\left(a\right)}},\,  \betaa{-1}{\overline{A}}{\rhoo{\betaa{}{\bar{a}}{\left(v\right)}}{\left(u\right)}} \right) \right),
\end{align*}
where
\begin{align*}
	\bar{a}:=\alphaa{-1}{u}{\left(a\right)}, \ &\bar{u}:= \betaa{-1}{a}{\left(u\right)},\\ A:=\alphaa{}{u}{\lambdaa{\bar{a}}{\left(b\right)}}, \ &U:=\betaa{}{a}{\lambdaa{\bar{u}}{\left(v\right)}}, \\ \overline{A}:=\alphaa{-1}{U}{\left(A\right)}, \ &\overline{U}:= \betaa{-1}{A}{\left(U\right)},
\end{align*}
for all $\left(a,u\right),\left(b,v\right)\in S\times T$, is a solution. This solution is called the \emph{matched product of the solutions} $r_S$ and $r_T$ (via $\alpha$ and $\beta$) and it is denoted by $r_S\bowtie r_T$.

\noindent If $\left(r_{S}, r_{T}, \alpha, \beta\right)$ is a matched product system of solutions, we denote $\alphaa{-1}{u}{\left(a\right)}$ with $\bar{a}$ and $\betaa{-1}{a}{\left(u\right)}$ with $\bar{u}$, when the pair $\left(a,u\right) \in S \times T$ is
clear from the context.

\section{The matched product of solutions associated with left shelves}

This section focuses on the matched product of solutions associated with left shelves. First, we prove that we can simplify the requirements to provide the matched product of solutions related to left shelves. In this case, the matched product solution is always left non-degenerate. Finally, we prove that the structure shelf of the matched product solution does not depend on the maps $\alpha$ and $\beta$.

\begin{theor}\label{th:2l-sh}
	Let $\left(S, \triangleright\right)$, $\left(T, \triangleright\right)$ be left shelves, $r_S$ and $r_T$ the solutions associated with $\left(S, \triangleright\right)$ and $\left(T, \triangleright\right)$, respectively, and $\alpha : T\to \Sym\left(S\right)$, $\beta : S \to \Sym\left(T\right)$ maps. Set $\alphaa{}{u}{}=\alpha\left(u\right)$ and $\betaa{}{a}{}=\betaa{}{}{\left(a\right)}$, for all $a\in S$ and $u \in T$. The quadruple $\left(r_{S}, r_{T}, \alpha, \beta\right)$ is a matched product system of solutions if and only if $\alphaa{}{u}{}$ and $\betaa{}{a}{}$ are homomorphisms of left shelves and 
	\begin{center}
		\begin{minipage}[b]{.5\textwidth}
			\vspace{-\baselineskip}
			\begin{align}\tag{l1}\label{eq:leftShelves1}
				\alpha_{v\triangleright u} =\alphaa{-1}{v}{\alphaa{}{u}{\alphaa{}{v}{}}}
			\end{align}
		\end{minipage}%
		\hfill\hfill\hfill
		\begin{minipage}[b]{.5\textwidth}
			\vspace{-\baselineskip}
			\begin{align}\tag{l2}\label{eq:leftShelves2}
				\beta_{b\triangleright a} = \betaa{-1}{b}{\betaa{}{a}{\betaa{}{b}{}}}
			\end{align}
		\end{minipage}
	\end{center}
	\begin{center}
		\begin{minipage}[b]{.5\textwidth}
			\vspace{-\baselineskip}
			\begin{align}\tag{l3}\label{eq:leftShelves3}
				\alphaa{}{u}{} = \alphaa{}{\betaa{-1}{a}{\left(u\right)}}{}
			\end{align}
		\end{minipage}%
		\hfill\hfill\hfill
		\begin{minipage}[b]{.5\textwidth}
			\vspace{-\baselineskip}
			\begin{align}\label{eq:leftShelves4}\tag{l4}
				\betaa{}{a}{} = \betaa{}{\alphaa{-1}{u}{\left(a\right)}}{}
			\end{align}
		\end{minipage}
	\end{center}
	hold for all $a,b \in S$, $u,v \in T$. In particular, the matched product solution $r_S\bowtie r_T$ is given by
	\begin{align*}
		r_S\bowtie r_T\left(\left(a,u\right),\left(b,v\right)\right) = \left(\left(\alphaa{}{u}{\left(b\right)},\betaa{}{a}{\left(v\right)}\right), 
		\left(\alphaa{-1}{v}{\left(\alphaa{}{u}{\left(b\right)}\triangleright a\right)},\betaa{-1}{b}{\left(\betaa{}{a}{\left(v\right)}\triangleright u\right)} \right)
		\right)
	\end{align*}
	and $r_S\bowtie r_T$ is left non-degenerate.
	
	\begin{proof}
		First, suppose that $\left(r_S,r_T,\alpha,\beta\right)$ is a matched product system of solutions. Since $r_S$ and $r_T$ are solutions associated with left shelves we have that $\lambdaa{a}{\left(b\right)} = b$, $\lambdaa{u}{\left(v\right)} = v$, $\rhoo{b}{\left(a\right)} = b\triangleright a$, and $\rhoo{v}{\left(u\right)} = v\triangleright u$  for all $a,b\in S$ and $u,v \in T$. Hence, by \eqref{eq:primo}, i.e., $\alphaa{}{u}{\alphaa{}{v}{}} = \alphaa{}{\lambdaa{u}{\left(v\right)}}{\alphaa{}{\rhoo{v}{\left(u\right)}}{}}$ we obtain $\alphaa{}{u}{\alphaa{}{v}{}}=\alphaa{}{v}{\alphaa{}{v\triangleright u}{}}$, i.e., \eqref{eq:leftShelves1} holds. Likewise, by \eqref{eq:secondo}, we obtain \eqref{eq:leftShelves2}. 
		Moreover, by \eqref{eq:terzo}, i.e., $\lambdaa{a}{\alphaa{}{\betaa{-1}{a}{\left(u\right)}}{}}=\alphaa{}{u}{\lambdaa{\alphaa{-1}{u}{\left(a\right)}}{}}$ we obtain \eqref{eq:leftShelves3} and, by \eqref{eq:quarto}, we have \eqref{eq:leftShelves4}.
		Furthermore, it follows that
		\begin{align*}
			\alphaa{-1}{u}{\left(b\right)}\triangleright\alphaa{-1}{u}{\left(a\right)}
			&=\alphaa{-1}{u}{\left(b\right)}\triangleright\alphaa{-1}{\betaa{}{a}{\left(u\right)}}{\left(a\right)}&\mbox{by \eqref{eq:leftShelves3}}\\
			&= \alpha^{-1}_{\beta_{\rho_b\left(a\right)}\beta^{-1}_b\left(u\right)}\rho_b\left(a\right)&\mbox{by \eqref{eq:quinto}}\\
			&=\alphaa{-1}{\betaa{}{b\triangleright a}{\betaa{-1}{b}{\left(u\right)}}}{\left(b\triangleright a\right)}\\
			&=\alphaa{-1}{\betaa{-1}{b}{\betaa{}{a}{\left(u\right)}}}{\left(b\triangleright a\right)}&\mbox{by \eqref{eq:leftShelves2}}\\
			&=\alphaa{-1}{u}{\left(b\triangleright a \right)}, &\mbox{by \eqref{eq:leftShelves3}}
		\end{align*}
		in other words, $b \triangleright a = \alphaa{-1}{u}{\left(\alphaa{}{u}{\left(b\right)}\triangleright\alphaa{}{u}{\left(a\right)} \right)}$.
		The previous equality proves that $\alphaa{}{u}{}$ is a homomorphism of left shelves. Likewise, one can see that $\betaa{}{a}{}$ is a homomorphism of left shelves. \\
		Conversely, if $\left(S,\triangleright\right)$ and $\left(T,\triangleright\right)$ are left shelves and $\alpha: T\to \Sym\left(S\right)$, $\beta:S\to \Sym\left(T\right)$ are maps that satisfy \eqref{eq:leftShelves1}, \eqref{eq:leftShelves2}, \eqref{eq:leftShelves3} and \eqref{eq:leftShelves4} and such that $\alpha_u$ and $\beta_a$ are left shelf homomorphisms for all $a\in S$ and $u\in T$, then $\left(r_S,r_T,\alpha,\beta\right)$ is a matched product system of solutions. Indeed, clearly \eqref{eq:leftShelves1}, \eqref{eq:leftShelves2}, \eqref{eq:leftShelves3} and \eqref{eq:leftShelves4} are equivalent to \eqref{eq:primo},\eqref{eq:secondo},\eqref{eq:terzo} and \eqref{eq:quarto}. Moreover, 
		\begin{align*}
			\alpha^{-1}_{\beta_{\rho_b\left(a\right)}\beta^{-1}_b\left(u\right)}\rho_b\left(a\right)
			&=\alphaa{-1}{\betaa{}{b\triangleright a}{\betaa{-1}{b}{\left(u\right)}}}{\left(b\triangleright a\right)}\\
			&=\alphaa{-1}{u}{\left(b\triangleright a\right)}&\mbox{by \eqref{eq:leftShelves3}}\\
			&=\alphaa{-1}{u}{\left(b\right)}\triangleright\alphaa{-1}{u}{\left(a\right)}&\mbox{since $\alpha_u$ is a homomorphism}\\
			&=\alphaa{-1}{u}{\left(b\right)}\triangleright\alphaa{-1}{\betaa{}{a}{\left(u\right)}}{\left(a\right)}&\mbox{by \eqref{eq:leftShelves3}}\\
			&=\rhoo{\alphaa{-1}{u}{\left(b\right)}}{\alphaa{-1}{\betaa{}{a}{\left(u\right)}}{\left(a\right)}},
		\end{align*}
		i.e., \eqref{eq:quinto} holds. By a similar argument, we also obtain that \eqref{eq:sesto} holds. Therefore, $\left(r_S,r_T,\alpha,\beta\right)$ is a matched product system of solutions.\\
		We compute $r_S\bowtie r_T$. First, straightforward, we have that 
		\begin{align*}
			\lambdaa{\left(a,u\right)}{\left(b,v\right)} 
			&=\left(\alphaa{}{u}{\lambdaa{\alphaa{-1}{u}{\left(a\right)}}{\left(b\right)}},\betaa{}{a}{\lambdaa{\betaa{-1}{a}{\left(u\right)}}{\left(v\right)}} \right)\\
			&=\left(\alphaa{}{u}{\left(b\right)},\betaa{}{a}{\left(v\right)}\right)
		\end{align*}
		Moreover, by setting $A:=\alphaa{}{u}{\left(b\right)}$ and $U:=\betaa{}{a}{\left(v\right)}$,
		we obtain that
		\begin{align*}
			\rhoo{\left(b,v\right)}{\left(a,u\right)} 
			&= \left(\alphaa{-1}{\overline{U}}{\rhoo{\alphaa{}{\bar{u}}{\left(b\right)}}{\left(a\right)}} , \betaa{-1}{\overline{A}}{\rhoo{\betaa{}{\bar{a}}{\left(v\right)}}{\left(u\right)}}\right)\\
			&=\left(\alphaa{-1}{\betaa{-1}{\alphaa{}{u}{\left(a\right)}}{\betaa{}{a}{\left(v\right)}}}{\rhoo{\alphaa{}{\betaa{-1}{a}{\left(u\right)}}{\left(b\right)}}{\left(a\right)}},
			\betaa{-1}{\alphaa{-1}{\betaa{}{a}{\left(u\right)}}{\alphaa{}{u}{\left(b\right)}}}{\rhoo{\betaa{}{\alphaa{-1}{u}{\left(a\right)}}{\left(v\right)}}{\left(u\right)}}\right)\\
			&=\left(\alphaa{-1}{\betaa{-1}{\alphaa{}{u}{\left(a\right)}}{\betaa{}{a}{\left(v\right)}}}{\left(\alphaa{}{\betaa{-1}{a}{\left(u\right)}}{\left(b\right)}\triangleright a\right)},
			\betaa{-1}{\alphaa{-1}{\betaa{}{a}{\left(u\right)}}{\alphaa{}{u}{\left(b\right)}}}{\left(\betaa{}{\alphaa{-1}{u}{\left(a\right)}}{\left(v\right)}\triangleright u\right)} \right)\\
			&=\left(\alphaa{-1}{v}{\left(\alphaa{}{u}{\left(b\right)}\triangleright a\right)},\betaa{-1}{b}{\left(\betaa{}{a}{\left(v\right)}\triangleright u\right)} \right),
		\end{align*}
		for all $a,b\in S$, $u,v \in T$.
		Finally, $\lambdaa{\left(a,u\right)}{}$ is bijective with inverse given by
		\begin{align*}
			\lambdaam{\left(a,u\right)}{\left(b,v\right)} =\left(\alphaa{-1}{u}{\left(b\right)},\betaa{-1}{a}{\left(v\right)}\right),
		\end{align*}
		i.e., $r_S\bowtie r_T$ is left non-degenerate.
	\end{proof}
\end{theor}

The following example shows that the matched product of two solutions associated with two left shelves is not necessarily associated with a left shelf.
\begin{ex}
	Let $\left(S, \triangleright\right)$ be the left shelf defined by $a\triangleright b:= a$, for all $a,b\in S$, and $\left(T, \triangleright\right) = \left(S, \triangleright \right)$.
	Then, the solution associated with $\left(S, \triangleright\right)$ is defined by $r_{S}\left(a, b\right) = \left(b, b\right)$, for all $a,b\in S$.
	Let $\theta$ and $\eta$ be a bijective maps from $S$ into itself, $\alpha,\beta:S\to \Sym\left(S\right)$ the constant maps with value $\theta$ and $\eta$, respectively. 
	Then $\left(r_{S}, r_{S}, \alpha, \beta\right)$ is a matched product system of solutions and $r_{S}\bowtie r_{S}$ is the map given by
	\begin{align*}
		r\left(\left(a,u\right), \left(b,v\right)\right)
		&= \left(\left(\theta\left(b\right),\eta\left(v\right)\right),\left( \theta^{-1}\left(\theta\left(b\right)\right), \eta^{-1}\left(\eta\left(v\right)\right)\right)\right)\\
		&=\left(\left(\theta\left(b\right),\eta\left(v\right)\right),\left(b,v\right)\right),
	\end{align*} 
	for all $a,b,u,v \in S$.
	Note that $r_{S}\bowtie r_{S}$ is associated with a left shelf if and only if $\lambda_{\left(a,u\right)}{\left(b,v\right)} = \left(b,v\right)$, for all $\left(a,u\right), \left(b,v\right)\in S\times T$ and, in other words, if and only if $\theta = \eta = \id_{S}$.
\end{ex}  

In the example above, we note that the solution associated with the matched product of two shelves is itself a solution related to a left shelf if the maps $\alpha$ and $\beta$ are trivial. This implication is true in general, as shown by the following theorem.

\begin{theor}\label{th:2l-sh->sh}
	Let $\left(S, \triangleright\right)$, $\left(T, \triangleright\right)$ be left shelves, $r_S$ and $r_T$ the solutions associated with $\left(S, \triangleright\right)$ and $\left(T, \triangleright\right)$, respectively. 
	If $\left(r_{S}, r_{T}, \alpha, \beta\right)$ is a matched product system of $r_{S}$ and $r_{T}$, then the matched solution $r_{S}\bowtie r_{T}$ is associated with a left shelf on the cartesian product $S\times T$ if and only if $\alpha_{u} = \id_{S}$, for every $u\in T$, and $\beta_{a} = \id_{T}$, for every $a\in S$.
	\begin{proof}
		First, suppose that the solution $r_{S}\bowtie r_{T}$ is associated with a left shelf on $S\times T$.
		It follows that
		\begin{align*}
			\left(\alphaa{}{u}{\left(b\right)},\betaa{}{a}{\left(u\right)}\right) 
			=\lambdaa{\left(a,u\right)}{\left(b,v\right)} =\left(b,v\right),
		\end{align*}
		for all $a,b\in S$, $u,v \in T$. Hence, 
		$\alphaa{}{u}{}=\id_S$ and $\betaa{}{a}{}=\id_{T}$
		for all $a \in S$, $u \in T$. 
		Conversely, suppose $\alpha_{u} = \id_{S}$, for every $u\in T$, and $\beta_{a} = \id_{T}$, for every $a\in S$. As a consequence,
		\begin{align*}
			\lambdaa{\left(a,u\right)}{\left(b,v\right)}
			=\left(\alphaa{}{u}{\left(b\right)},\betaa{}{a}{\left(v\right)}\right)
			=\left(b,v\right),
		\end{align*}
		for all $a,b\in S$, $u,v \in T$
		and, thus, the solution $r_{S}\bowtie r_{T}$ is a solution associated with a left shelf on $S\times T$.
	\end{proof}
\end{theor}

Finally, we compute the structure shelf of the matched product of two left shelves solutions, and we prove that the structure shelf does not depend on $\alpha$ and $\beta$.

\begin{prop}\label{prop:structureshelfleftleft}
	Let $\left(S,\triangleright\right)$, $\left(T,\triangleright\right)$ be left shelves, $r_S$ and $r_T$ the solutions associated with $\left(S,\triangleright\right)$ and  $\left(T,\triangleright\right)$, respectively. If $\left(r_{S}, r_{T}, \alpha, \beta\right)$ is a matched product system of $r_{S}$ and $r_{T}$, then the structure shelf of the matched solution $r:=r_{S}\bowtie r_{T}$ is  
	\begin{align*}
		\left(a,u\right)\triangleright_{r}\left(b,v\right) = \left(a\triangleright b, u\triangleright v\right),
	\end{align*}
	for all $a,b \in S$ and $u,v \in T$.
	\begin{proof}
		It follows from \cref{th:shleft}.
	\end{proof}
\end{prop}

\section{The matched product of solutions associated with right shelves}

In this section, we treat solutions associated with right shelves. We show that the requirements to obtain a matched product system of solutions are more straightforward. In this case, the matched product solution is left non-degenerate if and only if we start with right racks. The structure shelf of the matched product of two right racks does not depend on $\alpha$ and $\beta$.

\begin{theor}\label{theor:2r-sh}
	Let $\left(S, \triangleleft\right)$, $\left(T, \triangleleft\right)$ be right shelves, $r_S$ and $r_T$ the solutions associated with $\left(S, \triangleleft\right)$ and $\left(T, \triangleleft\right)$, respectively, and $\alpha : T\to \Sym\left(S\right)$, $\beta : S \to \Sym\left(T\right)$ maps. The quadruple $\left(r_{S}, r_{T}, \alpha, \beta\right)$ is a matched product system of solutions if and only if $\alphaa{}{u}{}$ and $\betaa{}{a}{}$ are homomorphisms of right shelves  and 
	\begin{center}
		\begin{minipage}[b]{.5\textwidth}
			\vspace{-\baselineskip}
			\begin{align}\tag{r1}\label{eq:rightShelves1}
				\alpha_{u\triangleleft v} = \alphaa{}{u}{\alphaa{}{v}{\alphaa{-1}{u}{}}}
			\end{align}
		\end{minipage}%
		\hfill\hfill\hfill
		\begin{minipage}[b]{.5\textwidth}
			\vspace{-\baselineskip}
			\begin{align}\tag{r2}\label{eq:rightShelves2}
				\beta_{a\triangleleft b} = \betaa{}{a}{\betaa{}{b}{\betaa{-1}{a}{}}}
			\end{align}
		\end{minipage}
	\end{center}
	\begin{center}
		\begin{minipage}[b]{.5\textwidth}
			\vspace{-\baselineskip}
			\begin{align}\tag{r3}\label{eq:rightShelves3}
				\alphaa{}{u}{} = \alphaa{}{\betaa{-1}{a}{\left(u\right)}}{}
			\end{align}
		\end{minipage}%
		\hfill\hfill\hfill
		\begin{minipage}[b]{.5\textwidth}
			\vspace{-\baselineskip}
			\begin{align}\label{eq:rightShelves4}\tag{r4}
				\betaa{}{a}{} = \betaa{}{\alphaa{-1}{u}{\left(a\right)}}{}
			\end{align}
		\end{minipage}
	\end{center}
	hold for all $a,b \in S$, $u,v \in T$. In particular, the matched product $r_S\bowtie r_T$ is given by
	\begin{align*}
		r_S\bowtie r_T\left(\left(a,u\right),\left(b,v\right)\right) = \left(\left(\alphaa{}{u}{\left(b\right)}\triangleleft a,\betaa{}{a}{\left(v\right)}\triangleleft u \right), \left(\alphaa{-1}{v\triangleleft u}{\left(a\right)},\betaa{-1}{b\triangleleft a}{\left(u\right)} \right) \right),
	\end{align*}
	for all $a,b \in S$, $u,v\in T$ and is left non-degenerate if and only if $\left(S,\triangleleft\right)$ and $\left(T,\triangleleft\right)$ are right racks.
	\begin{proof}
		First, suppose that $\left(r_S,r_T,\alpha,\beta\right)$ is a matched product system of solutions. Since $r_S$ and $r_T$ are solutions associated with right shelves we have that $\lambdaa{a}{\left(b\right)} = b\triangleleft a$, $\lambdaa{u}{\left(v\right)} = v\triangleleft u$, $\rhoo{b}{\left(a\right)} = a$, and $\rhoo{v}{\left(u\right)} = u $, for all $a,b\in S$ and $u,v \in T$. Hence,
		by \eqref{eq:primo} and \eqref{eq:secondo}, we obtain \eqref{eq:rightShelves1} and \eqref{eq:rightShelves2}. 
		Moreover, by \eqref{eq:quinto} we obtain
		\begin{align}\label{eq:rightShelfstar}
			\alphaa{-1}{\betaa{}{a}{\left(u\right)}}{\left(a\right)}
			= \alphaa{-1}{\betaa{}{a}{\betaa{-1}{b}{\left(u\right)}}}{\left(a\right)}
		\end{align}
		and 
		\begin{align}\label{eq:rightShelf2star}
			\alphaa{-1}{\betaa{}{a}{\left(u\right)}}{\left(a\right)} = \alphaa{-1}{u}{\left(a\right)}.
		\end{align}
		Hence, it follows that
		\begin{align*}
			\alphaa{-1}{u}{\left(b\right)} &=\alphaa{-1}{\betaa{}{b}{\left(u\right)}}{\left(b\right)} &\mbox{by \eqref{eq:rightShelf2star}}\\
			&=\alphaa{-1}{\betaa{}{b}{\betaa{-1}{a}{\left(u\right)}}}{\left(b\right)} &\mbox{by \eqref{eq:rightShelfstar}}\\
			&=\alphaa{-1}{\betaa{-1}{a}{\left(u\right)}}{\left(b\right)}, &\mbox{by \eqref{eq:rightShelf2star}}
		\end{align*}
		i.e., \eqref{eq:rightShelves3} holds. Likewise, we obtain that \eqref{eq:rightShelves4} is satisfied. Furthermore,
		\begin{align*}
			\alphaa{}{u}{\left(b\triangleleft \alphaa{-1}{u}{\left(a\right)} \right)} &=
			\alphaa{}{\betaa{-1}{a}{\left(u\right)}}{\left(b\right)}\triangleleft a &\mbox{by \eqref{eq:terzo}}\\
			&= \alphaa{}{u}{\left(b\right)}\triangleleft a. &\mbox{by \eqref{eq:rightShelves3}}
		\end{align*}
		It follows that $\alphaa{}{u}{\left(b\triangleleft a\right)}=\alphaa{}{u}{\left(b\right)}\triangleleft\alphaa{}{u}{\left(a\right)}$, i.e., $\alphaa{}{u}{}$ is a right shelf homomorphism. By \eqref{eq:quarto}, \eqref{eq:rightShelves4}, and a similar argument we prove that $\beta_a$ is a right shelf homomorphism.\\
		Conversely, suppose that $\left(S,\triangleleft\right)$ and $\left(T,\triangleleft\right)$ are right shelves and $\alpha: T\to \Sym\left(S\right)$, $\beta:S\to \Sym\left(T\right)$ are maps that satisfy \eqref{eq:rightShelves1}, \eqref{eq:rightShelves2}, \eqref{eq:rightShelves3} and \eqref{eq:rightShelves4} and such that $\alpha_u$ and $\beta_a$ are right shelf homomorphisms for all $a\in S$ and $u\in T$. Clearly, \eqref{eq:primo} and \eqref{eq:secondo} hold by \eqref{eq:rightShelves1} and \eqref{eq:rightShelves2}. Moreover,
		\begin{align*}
			\alphaa{-1}{\betaa{}{\rhoo{b}{\left(a\right)}}{\betaa{-1}{b}{\left(u\right)}}}{\rhoo{b}{\left(a\right)}}
			&=\alphaa{-1}{\betaa{}{a}{\betaa{-1}{b}{\left(u\right)}}}{\left(a\right)}\\
			&=\alphaa{-1}{\betaa{-1}{b}{\left(u\right)}}{\left(a\right)}&\mbox{by \eqref{eq:rightShelves3}}\\
			&=\alphaa{-1}{u}{\left(a\right)} &\mbox{by \eqref{eq:rightShelves3}}\\
			&=\alphaa{-1}{\betaa{}{a}{\left(u\right)}}{\left(a\right)} &\mbox{by \eqref{eq:rightShelves3}}\\
			&= \rhoo{\alphaa{-1}{u}{\left(b\right)}}{\alphaa{-1}{\betaa{}{a}{\left(u\right)}}{\left(a\right)}},
		\end{align*}
		i.e., \eqref{eq:quinto} holds. And with same computation \eqref{eq:sesto} holds. Furthermore, 
		\begin{align*}
			\lambdaa{a}{\alphaa{}{\betaa{-1}{a}{\left(u\right)}}{\left(b\right)}}
			&= \alphaa{}{\betaa{-1}{a}{\left(u\right)}}{\left(b\right)}\triangleleft a\\
			&=\alphaa{}{u}{\left(b\right)}\triangleleft a &\mbox{by \eqref{eq:rightShelves3}}\\
			&=\alphaa{}{u}{\left(b\triangleleft\alphaa{-1}{u}{\left(a\right)} \right)} &\mbox{since $\alpha_u$ is a homomorphism}\\
			&=\alphaa{}{u}{\lambdaa{\alphaa{-1}{u}{\left(a\right)}}{\left(b\right)}},
		\end{align*}
		i.e., \eqref{eq:terzo} holds. Similarly, \eqref{eq:quarto} is satisfied. Therefore, $\left(r_S,r_T,\alpha,\beta\right)$ is a matched product system of solutions. We compute the solution $r_S\bowtie r_T$:
		\begin{align*}
			\lambdaa{\left(a,u\right)}{\left(b,v\right)}
			&=\left(\alphaa{}{u}{\left(b\triangleleft\alphaa{-1}{u}{\left(a\right)}\right)},\betaa{}{a}{\left(v\triangleleft\betaa{-1}{a}{\left(u\right)}\right)}\right)\\
			&=\left(\alphaa{}{u}{\left(b\right)}\triangleleft a, \betaa{}{a}{\left(v\right)}\triangleleft u\right)
		\end{align*}
		and
		\begin{align*}
			\rhoo{\left(b,v\right)}{\left(a,u\right)} 
			&=\left(\alphaa{-1}{\betaa{-1}{\alphaa{}{u}{\left(b\right)}\triangleleft a}{\left(\betaa{}{a}{\left(v\right)}\triangleleft u\right)}}{\left(a\right)},\betaa{-1}{\alphaa{-1}{\betaa{}{a}{\left(v\right)}\triangleleft u}{\left(\alphaa{}{u}{\left(b\right)}\triangleleft a\right)}}{\left(u\right)} \right) \\
			&=\left(\alphaa{-1}{\betaa{}{a}{\left(v\right)\triangleleft u}}{\left(a\right)},\betaa{-1}{\alphaa{}{u}{\left(b\right)\triangleleft a}}{\left(u\right)} \right) &\mbox{by \eqref{eq:rightShelves3} and \eqref{eq:rightShelves4}}\\
			&=\left(\alphaa{}{\betaa{}{a}{\left(v\right)}}{\alphaa{-1}{u}{\alphaa{-1}{\betaa{}{a}{\left(v\right)}}{\left(a\right)}}},
			\betaa{}{\alphaa{}{u}{\left(b\right)}}{\betaa{-1}{a}{\betaa{-1}{\alphaa{}{u}{\left(b\right)}}{\left(u\right)}}}\right)&\mbox{by \eqref{eq:rightShelves1} and \eqref{eq:rightShelves2}}\\
			&=\left(\alphaa{}{v}{\alphaa{-1}{u}{\alphaa{-1}{v}{\left(a\right)}}}, \betaa{}{b}{\betaa{-1}{a}{\betaa{-1}{b}{\left(u\right)}}}\right)&\mbox{by \eqref{eq:rightShelves3} and \eqref{eq:rightShelves4}}\\
			&=\left(\alphaa{-1}{v\triangleleft u}{\left(a\right)}, \betaa{-1}{b\triangleleft a}{\left(u\right)}\right).&\mbox{by \eqref{eq:rightShelves1} and \eqref{eq:rightShelves2}}
		\end{align*}
		Finally, define $R_a\left(b\right):=b\triangleleft a$ and $R_u\left(v\right):= u\triangleleft v$.
		Suppose that $\left(S,\triangleleft\right)$ and $\left(T,\triangleleft\right)$ are right racks, then $R_a$ and $R_u$ are bijective and also $\lambdaa{\left(a,u\right)}{}$ is invertible with inverse given by
		\begin{align*}
			\lambdaam{\left(a,u\right)}{\left(b,v\right)}=\left(\alphaa{-1}{u}{R^{-1}_a\left(b\right)},\betaa{-1}{a}{R^{-1}_u\left(v\right)}\right).
		\end{align*}
		Conversely, if $\lambdaa{\left(a,u\right)}{}$ is bijective then for a fixed pair $\left(c,w\right)\in S\times T$ there exists a pair $\left(c,w\right)$ such that $\left(c,w\right) = \lambdaa{\left(a,u\right)}{\left(b,w\right)} = \left(R_a\alphaa{}{u}{\left(b\right)}, R_u\betaa{}{a}{\left(v\right)}\right)$ and so $c = R_a\alphaa{}{u}{\left(b\right)}$ and $w = R_u\betaa{}{a}{\left(v\right)}$ and $R_a$ and $R_u$ are surjective. Now, suppose that $b,c\in S$ and $v,w\in T$ such that $R_a\left(b\right) = R_a\left(c\right)$ and $R_u\left(v\right)=R_u\left(w\right)$, then
		\begin{align*}
			\lambdaa{\left(a,u\right)}{\left(\alphaa{-1}{u}{\left(b\right)}, \betaa{-1}{a}{\left(v\right)}\right)}
			&=\left(R_a\left(b\right),R_u\left(v\right)\right)\\
			&=\left(R_a\left(c\right),R_u\left(w\right)\right)\\
			&=\lambdaa{\left(a,u\right)}{\left(\alphaa{-1}{u}{\left(c\right)}, \betaa{-1}{a}{\left(w\right)}\right)}.
		\end{align*}
		Since $\lambdaa{\left(a,u\right)}{}$ is injective then $\alphaa{-1}{u}{\left(b\right)} = \alphaa{-1}{u}{\left(c\right)}$ and $\betaa{-1}{a}{\left(v\right)} =\betaa{-1}{a}{\left(w\right)}$ and so $b=c$ and $v=w$. Hence, $R_a$ and $R_u$ are also injective and $\left(S,\triangleleft\right)$ and $\left(T,\triangleright\right)$ are right racks.
		This completes the proof.
	\end{proof}
\end{theor}

Moreover, if $\left(S,\triangleleft\right)$, $\left(T,\triangleleft\right)$ are right shelves and $\alpha:T\to \Sym\left(S\right)$, $\beta:S\to \Sym\left(T\right)$ are maps such that $\left(r_S,r_T,\alpha,\beta\right)$ is a matched product system of the solutions $r_S$, $r_T$ associated with $\left(S,\triangleleft\right)$, $\left(T,\triangleleft\right)$ respectively, then $r_S\bowtie r_T$ is associated with a right shelf if and only if $\rhoo{\left(b,v\right)}{\left(a,u\right)} = \left(a,u\right)$. By \cref{theor:2r-sh}, this means that $\alphaa{-1}{v\triangleright u}{\left(a\right)}=a$ and $\betaa{}{b\triangleright a}{\left(u\right)}=u$, for all $a,b\in S$, $u,v \in T$. By \eqref{eq:rightShelves1} and \eqref{eq:rightShelves2} we obtain that $\alphaa{-1}{v}{\alphaa{-1}{u}{\left(a\right)}}=\alphaa{-1}{u}{\left(a\right)}$ and $\betaa{-1}{b}{\betaa{-1}{a}{\left(u\right)}} = \betaa{-1}{a}{\left(u\right)}$, i.e., the following theorem holds.
\begin{theor}\label{th:2r-sh->sh}
	Let $\left(S, \triangleleft\right)$, $\left(T, \triangleleft\right)$ be right shelves, $r_S$ and $r_T$ the solutions associated with $\left(S, \triangleright\right)$ and $\left(T, \triangleright\right)$, respectively. 
	If $\left(r_{S}, r_{T}, \alpha, \beta\right)$ is a matched product system of $r_{S}$ and $r_{T}$, then the matched product solution $r_{S}\bowtie r_{T}$ is associated with a right shelf on the cartesian product $S\times T$ if and only if $\alpha_{u} = \id_{S}$, for every $u\in T$, and $\beta_{a} = \id_{T}$, for every $a\in S$.
\end{theor}

\begin{prop}\label{prop:structureshelfrightright}
	Let $\left(S,\triangleleft\right)$, $\left(T,\triangleleft\right)$ be right racks, $r_S$ and $r_T$ the solutions associated with $\left(S,\triangleleft\right)$ and  $\left(T,\triangleleft\right)$, respectively. If $\left(r_{S}, r_{T}, \alpha, \beta\right)$ is a matched product system of $r_{S}$ and $r_{T}$, then the structure shelf of the matched solution $r:=r_{S}\bowtie r_{T}$ is 
	\begin{align*}
		\left(a,u\right)\triangleright_{r} \left(b,v\right) := \left(b\triangleleft a, v \triangleleft u\right),
	\end{align*}
	for all $a,b \in S$ and $u,v \in T$.
	\begin{proof}
		It follows from \cref{th:shleft}.
	\end{proof}
\end{prop}

\section{The matched product of a solution associated with a left shelf and a solution associated with a right shelf}

In this section, we present the matched product of one solution associated with a left shelf and the other associated with a right shelf. In this case, we also obtain natural conditions to check the properties of the matched product system. Then, we give conditions for left non-degeneracy and, under these assumptions, we provide the structure shelf of the matched product solution.

\begin{theor}\label{theor:lr-sh}
	Let $\left(S, \triangleright\right)$ be a left shelf, $\left(T, \triangleleft\right)$ be right shelf, $r_S$ and $r_T$ the solutions associated with $\left(S, \triangleright\right)$ and $\left(T, \triangleleft\right)$, respectively, and $\alpha : T\to \Sym\left(S\right)$, $\beta : S \to \Sym\left(T\right)$ maps. The quadruple $\left(r_{S}, r_{T}, \alpha, \beta\right)$ is a matched product system of solutions if and only if $\alphaa{}{u}{}$ are homomorphisms of left shelves and $\betaa{}{a}{}$ are homomorphisms of right shelves for all $a\in S$, $u \in T$ and 
	\begin{center}
		\begin{minipage}[b]{.5\textwidth}
			\vspace{-\baselineskip}
			\begin{align}\tag{lr1}\label{eq:leftrightShelves1}
				\alpha_{u\triangleleft v} = \alpha_{u}\alpha_{v}\alpha_{u}^{-1}
			\end{align}
		\end{minipage}%
		\hfill\hfill\hfill
		\begin{minipage}[b]{.5\textwidth}
			\vspace{-\baselineskip}
			\begin{align}\tag{lr2}\label{eq:leftrightShelves2}
				\beta_{a\triangleright b} = \betaa{-1}{b}{\betaa{}{a}{\betaa{}{b}{}}}
			\end{align}
		\end{minipage}
	\end{center}
	\begin{center}
		\begin{minipage}[b]{.5\textwidth}
			\vspace{-\baselineskip}
			\begin{align}\tag{lr3}\label{eq:leftrightShelves3}
				\alphaa{}{u}{} = \alphaa{}{\betaa{-1}{a}{\left(u\right)}}{}
			\end{align}
		\end{minipage}%
		\hfill\hfill\hfill
		\begin{minipage}[b]{.5\textwidth}
			\vspace{-\baselineskip}
			\begin{align}\label{eq:leftrightShelves4}\tag{lr4}
				\betaa{}{a}{} = \betaa{}{\alphaa{-1}{u}{\left(a\right)}}{}
			\end{align}
		\end{minipage}
	\end{center}
	hold for all $a,b \in S$, $u,v \in T$. In particular, the matched product $r_S\bowtie r_T$ is given by
	\begin{align*}
		r_S\bowtie r_T\left(\left(a,u\right),\left(b,v\right)\right) 
		= \left(
		\left(\alphaa{}{u}{\left(b\right)}, \betaa{}{a}{\left(v\right)}\triangleleft\ u\right), 
		\left(\alphaa{-1}{v\triangleleft u}{\alphaa{}{u}{\left(b \right)}\triangleright a},  \betaa{-1}{b}{\left(u \right)} \right) 
		\right),
	\end{align*}
	for all $a,b \in S$, $u,v \in T$ and $r_S\bowtie r_T$ is left non-degenerate if and only if $\left(T,\triangleleft\right)$ is a right rack.
	\begin{proof}
		First, suppose that $\left(r_S,r_T,\alpha,\beta\right)$ is a matched product system of solutions. Since $r_S$ is a solution associated with a left shelf and $r_T$ is associated with a right shelf we have that $\lambdaa{a}{\left(b\right)} =b$, $\lambdaa{u}{\left(v\right)} = v\triangleleft u$, $\rhoo{b}{\left(a\right)} = b\triangleright a$, and $\rhoo{v}{\left(u\right)} = u $  for all $a,b\in S$ and $u,v \in T$. Hence,
		by \eqref{eq:primo} and \eqref{eq:secondo} we obtain \eqref{eq:leftrightShelves1} and \eqref{eq:leftrightShelves2}.			
		Moreover, by \eqref{eq:terzo} we have that 
		\begin{align*}
			\alphaa{}{\betaa{-1}{a}{\left(u\right)}}{\left(b\right)} = \alphaa{}{u}{\left(b\right)},
		\end{align*}
		i.e., \eqref{eq:leftrightShelves3} holds. 
		Furthermore, we obtain that
		\begin{align*}
			\alphaa{-1}{u}{\left(b\right)}\triangleright\alphaa{-1}{u}{\left(a\right)}&=\alphaa{-1}{u}{\left(b\right)}\triangleright \alphaa{-1}{\betaa{}{a}{\left(u\right)}}{\left(a\right)} &\mbox{by \eqref{eq:leftrightShelves3}}\\
			&=\alphaa{-1}{\betaa{}{b\triangleright a}{\betaa{-1}{b}{\left(u\right)}}}{\left(b\triangleright a\right)}&\mbox{by \eqref{eq:quinto}}\\
			&=\alphaa{-1}{u}{\left(b\triangleright a\right)}&\mbox{by \eqref{eq:leftrightShelves3}}
		\end{align*}
		and then $\alphaa{}{u}{}$ is a right shelf homomorphism. Moreover, by \eqref{eq:sesto} we have that
		\begin{align}\label{eq:alphabeta3}
			\betaa{-1}{\alphaa{}{u}{\left(a\right)}}{\left(u\right)} = \betaa{-1}{\alphaa{}{u}{\alphaa{-1}{v}{\left(a\right)}}}{\left(u\right)},
		\end{align} 
		for all $a \in S$, $u,v \in T$. Therefore,
		\begin{align}\label{eq:alphabeta4}
			\betaa{-1}{a}{\left(u\right)} = \betaa{-1}{\alphaa{}{u}{\alphaa{-1}{u}{\left(a\right)}}}{\left(u\right)} = \betaa{-1}{\alphaa{}{u}{\left(a\right)}}{\left(u\right)},
		\end{align}
		for all $a \in S$ and $u \in T$. It follows that 
		\begin{align*}
			\betaa{-1}{a}{\left(v\right)} &= \betaa{-1}{\alphaa{}{v}{\left(a\right)}}{\left(v\right)} &\mbox{by \eqref{eq:alphabeta4}}\\
			&= \betaa{-1}{\alphaa{}{v}{\alphaa{-1}{u}{\left(a\right)}}}{\left(v\right)} &\mbox{by \eqref{eq:alphabeta3}}\\
			&= \betaa{-1}{\alphaa{-1}{u}{\left(a\right)}}{\left(v\right)},&\mbox{by \eqref{eq:alphabeta4}}
		\end{align*}
		for all $a \in S$ and $u,v \in T$, i.e., \eqref{eq:leftrightShelves4} holds. 
		Furthermore, for all $a\in S$ and $u,v \in T$ we obtain that
		\begin{align*}
			\betaa{}{a}{\left(v\triangleright\betaa{-1}{a}{\left(u\right)}\right)}&=\betaa{}{\alphaa{-1}{u}{\left(a\right)}}{\left(v\right)} \triangleright u &\mbox{by \eqref{eq:quarto}}\\
			&=\betaa{}{a}{\left(v\right)} \triangleright u &\mbox{by \eqref{eq:leftrightShelves4}}
		\end{align*}
		and then $\betaa{}{a}{\left(v\triangleright u\right)} = \betaa{}{a}{\left(v\right)}\triangleright\betaa{}{a}{\left(u\right)}$, for all $a \in S$, $u,v\in T$, i.e., $\betaa{}{a}{}$ is a right shelf homomorphism, for every $a \in S$.
		Therefore, if $\left(r_S,r_T,\alpha,\beta\right)$ is a matched product system of solutions then $\alphaa{}{u}{}$ is a left shelf homomorphism, for every $u\in T$, $\betaa{}{a}{}$ is a right shelf homomorphism, for every $a\in S$, and the properties \eqref{eq:leftrightShelves1}, \eqref{eq:leftrightShelves2}, \eqref{eq:leftrightShelves3}, and \eqref{eq:leftrightShelves4} hold.\\
		Conversely, suppose that $\left(S,\triangleright\right)$ is a left shelf, $\left(T,\triangleleft\right)$ is a right shelf, $\alpha: T\to \Sym\left(S\right)$, $\betaa{}{}{}:S\to \Sym\left(T\right)$ are such that $\alphaa{}{u}{}$ is a left shelf homomorphism for every $u\in T$, $\betaa{}{a}{}$ is a right shelf homomorphism for every $a\in S$, and properties \eqref{eq:leftrightShelves1}, \eqref{eq:leftrightShelves2}, \eqref{eq:leftrightShelves3}, and \eqref{eq:leftrightShelves4} hold. Clearly, \eqref{eq:primo} and \eqref{eq:secondo} are satisfied. Furthermore, 
		\begin{align*}
			\rhoo{\alphaa{-1}{u}{\left(b\right)}}{\alphaa{-1}{\betaa{}{a}{\left(u\right)}}{\left(a\right)}}
			&=\alphaa{-1}{u}{\left(b\right)}\triangleright \alphaa{-1}{\betaa{}{a}{\left(u\right)}}{\left(a\right)}\\
			&=\alphaa{-1}{u}{\left(b\right)}\triangleright\alphaa{-1}{u}{\left(a\right)}&\mbox{by \eqref{eq:leftrightShelves3}}\\
			&=\alphaa{-1}{u}{\left(b\triangleright a\right)} &\mbox{$\alpha_u$ is a homomorphism}\\
			&=\alphaa{-1}{\betaa{}{b\triangleright a}{\betaa{-1}{b}{\left(u\right)}}}{\left(b \triangleright a\right)}&\mbox{by \eqref{eq:leftrightShelves3}}\\
			&=\alphaa{-1}{\betaa{}{\rhoo{b}{\left(a\right)}}{\betaa{-1}{b}{\left(u\right)}}}{\rhoo{b}{\left(a\right)}},
		\end{align*}
		i.e., \eqref{eq:quinto} holds. In addition,
		\begin{align*}
			\rhoo{\betaa{-1}{a}{\left(v\right)}}{\betaa{-1}{\alphaa{}{u}{\left(a\right)}}{\left(u\right)}}
			&=\betaa{-1}{\alphaa{}{u}{\left(a\right)}}{\left(u\right)}\\
			&=\betaa{-1}{a}{\left(u\right)} &\mbox{by \eqref{eq:leftrightShelves4}}\\
			&=\betaa{-1}{\alphaa{}{u}{\alphaa{-1}{v}{\left(a\right)}}}{\left(u\right)} &\mbox{by \eqref{eq:leftrightShelves4}}\\
			&=\betaa{-1}{\alphaa{}{\rhoo{v}{\left(u\right)}}{\alphaa{-1}{v}{\left(a\right)}}}{\left(u\right)},
		\end{align*}
		i.e., \eqref{eq:sesto} is satisfied. Next,
		\begin{align*}
			\lambdaa{a}{\alphaa{}{\betaa{-1}{a}{\left(u\right)}}{\left(b\right)}}
			&=\alphaa{}{\betaa{-1}{a}{\left(u\right)}}{\left(b\right)} \\
			&=\alphaa{}{u}{\left(b\right)}&\mbox{by \eqref{eq:leftrightShelves3}}\\
			&=\alphaa{}{u}{\lambdaa{\alphaa{-1}{u}{\left(a\right)}}{\left(b\right)}}
		\end{align*}
		and
		\begin{align*}
			\lambdaa{u}{\betaa{}{\alphaa{-1}{u}{\left(a\right)}}{\left(v\right)}}
			&=\betaa{}{\alphaa{-1}{u}{\left(a\right)}}{\left(v\right)}\triangleleft u\\
			&=\betaa{}{a}{\left(v\right)}\triangleleft u &\mbox{by \eqref{eq:leftrightShelves4}}\\
			&=\betaa{}{a}{\left(v \triangleleft \betaa{-1}{a}{\left(u\right)}\right)} &\mbox{$\betaa{}{a}{}$ in a homomorphism}\\
			&=\betaa{}{a}{\lambdaa{\betaa{-1}{a}{\left(u\right)}}{\left(v\right)}}.
		\end{align*}
		Hence, \eqref{eq:terzo} and \eqref{eq:quarto} hold. It follows that $\left(r_S,r_T,\alpha, \beta\right)$ is a matched product system of solutions.\\
		Furthermore, we can compute the matched product solution, $r_S\bowtie r_T$. It is given by $r_S\bowtie r_T\left(\left(a,u\right),\left(b,v\right)\right)=\left(\lambdaa{\left(a,u\right)}{\left(b,v\right)},\rhoo{\left(b,v\right)}{\left(a,u\right)}\right)$ where
		\begin{align*}
			\lambdaa{\left(a,u\right)}{\left(b,v\right)}
			&=\left(\alphaa{}{u}{\lambdaa{\alphaa{-1}{u}{\left(a\right)}}{\left(b\right)}}, \betaa{}{a}{\lambdaa{\betaa{-1}{a}{\left(u\right)}}{\left(v\right)}}\right)\\
			&=\left(\alphaa{}{u}{\left(b\right)},\betaa{}{a}{\left(v\triangleleft\betaa{-1}{a}{\left(u\right)}\right)}\right)\\
			&=\left(\alphaa{}{u}{\left(b\right)},\betaa{}{a}{\left(v\right)}\triangleleft u\right)
		\end{align*}
		and, by setting $A:=\alphaa{}{u}{\left(b\right)}$ and $U:=\betaa{}{a}{\left(v\right)}\triangleleft u$,
		\begin{align*}
			\rhoo{\left(b,v\right)}{\left(a,u\right)} 
			&=\left(\alphaa{-1}{\betaa{-1}{A}{\left(U\right)}}{\rhoo{\alphaa{}{\betaa{-1}{a}{\left(u\right)}}{\left(b\right)}}{\left(a\right)}}, \betaa{-1}{\alphaa{-1}{U}{\left(A\right)}}{\rhoo{\betaa{}{\alphaa{-1}{u}{\left(a\right)}}{\left(v\right)}}{\left(u\right)}}\right)\\
			&= \left(\alphaa{-1}{U}{\rhoo{\alphaa{}{u}{\left(b\right)}}{\left(a\right)}}, \betaa{-1}{A}{\rhoo{\betaa{}{a}{\left(v\right)}}{\left(u\right)}}\right) &\mbox{by \eqref{eq:leftrightShelves4} and \eqref{eq:leftrightShelves3}}\\
			&=\left(\alphaa{-1}{\betaa{}{a}{\left(v\right)}\triangleleft u}{\left(\alphaa{}{u}{\left(b\right)}\triangleright a\right)},\betaa{-1}{\alphaa{}{u}{\left(b\right)}}{\left(u\right)} \right)\\
			&=\left(\alphaa{}{\betaa{}{a}{\left(v\right)}}{\alphaa{-1}{u}{\alphaa{-1}{\betaa{}{a}{\left(v\right)}}{\left(\alphaa{}{u}{\left(b\right)}\triangleright a\right)}}}, \betaa{-1}{b}{\left(u\right)} \right)&\mbox{by \eqref{eq:leftrightShelves1} and \eqref{eq:leftrightShelves4}}\\
			&=\left(\alphaa{}{v}{\alphaa{-1}{u}{\alphaa{-1}{v}{\left(\alphaa{}{u}{\left(b\right)}\triangleright a\right)}}},\betaa{-1}{b}{\left(u\right)} \right)&\mbox{by \eqref{eq:leftrightShelves3}}\\
			&=\left(\alphaa{-1}{v\triangleleft u}{\alphaa{}{u}{\left(b\right)\triangleright a}}, \betaa{-1}{b}{\left(u\right)} \right)&\mbox{by \eqref{eq:leftrightShelves1}}.
		\end{align*}
		Finally, it is easy to check that $r_S\bowtie r_T$ is left non-degenerate, i.e., $\lambdaa{\left(a,u\right)}{}$ is bijective if and only if $\left(T,\triangleleft\right)$ is a right rack.
	\end{proof}
\end{theor}

\begin{prop}\label{prop:structureshelfleftright}
	Let $\left(S,\triangleright\right)$ be a left shelf, $\left(T,\triangleleft\right)$ be right rack, $r_S$ and $r_T$ the solutions associated with $\left(S,\triangleright\right)$ and  $\left(T,\triangleleft\right)$, respectively. If $\left(r_{S}, r_{T}, \alpha, \beta\right)$ is a matched product system of $r_{S}$ and $r_{T}$, then the structure shelf of the matched solution $r:=r_{S}\bowtie r_{T}$ is
	\begin{align*}
		\left(a,u\right)\triangleright_{r} \left(b,v\right) := \left(a\triangleright b, v \triangleleft u\right),
	\end{align*}
	for all $a,b \in S$ and $u,v \in T$.
	\begin{proof}
		It follows from \cref{th:shleft}.
	\end{proof}
\end{prop}

\section{The structure shelf of the matched product of left non-degenerate solutions}

In this last section, we prove that the Propositions \ref{prop:structureshelfleftleft},   \ref{prop:structureshelfrightright} and \ref{prop:structureshelfleftright} are particular cases of a general result. Indeed, in the following theorem, we show that the structure shelf of the matched product of left non-degenerate solutions does not depend on $\alpha$ and $\beta$. 

\begin{theor}\label{th:shleft}
	Let $\left(r_S,r_T, \alpha,\beta\right)$ be a matched product system of left non-degenerate solutions. Then the matched product solution $r:=r_S\bowtie r_T$ is left non-degenerate and the structure shelf is 
	\begin{align*}
		\left(a,u\right)\triangleright_{r}\left(b,v\right)
		=\left(a\triangleright_{r_S}b, u\triangleright_{r_T}v \right),
	\end{align*}
	for all $a,b \in S$ and $u,v \in T$.
	\begin{proof}
		Recall that 
		\begin{align*}
			\lambdaa{\left(a,u\right)}{\left(b,v\right)} = \left(\alphaa{}{u}{\lambdaa{\bar{a}}{\left(b\right)}},\betaa{}{a}{\lambdaa{\bar{u}}{\left(v\right)}}\right)
		\end{align*}
		and
		\begin{align*}
			\rhoo{\left(b,v\right)}{\left(a,u\right)} = \left(\alphaa{-1}{\bar{U}}{\rhoo{\alphaa{}{\bar{u}}{\left(b\right)}}{\left(a\right)}},\betaa{-1}{\bar{A}}{\rhoo{\betaa{}{\bar{a}}{\left(v\right)}}{\left(u\right)}}\right),
		\end{align*}
		where $A = \alphaa{}{u}{\lambdaa{\bar{a}}{\left(b\right)}}$ and $U=\betaa{}{a}{\lambdaa{\bar{u}}{\left(v\right)}}$.
		Since $r_S$ and $r_T$ are left non-degenerate, then $\lambdaa{\left(a,u\right)}{}$ is bijective with inverse given by
		\begin{align*}
			\lambdaam{\left(a,u\right)}{\left(b,v\right)}=\left(\lambdaam{\bar{a}}{\alphaa{-1}{u}{\left(b\right)}}, \lambdaam{\bar{u}}{\betaa{-1}{a}{\left(v\right)}}\right),
		\end{align*}
		for all $a,b\in S$ and $u,v\in T$. Hence we can compute the associated structure shelf, i.e.,
		\begin{align*}
			\left(a,u\right)&\triangleright_{r}\left(b,v\right) 
			=\lambdaa{\left(a,u\right)}{\rhoo{\lambdaam{\left(b,v\right)}{\left(a,u\right)}}{\left(b,v\right)}} \\
			&=\lambdaa{\left(a,u\right)}{\rhoo{\left(\lambdaam{\bar{b}}{\alphaa{-1}{v}{\left(a\right)}},\lambdaam{\bar{v}}{\betaa{-1}{b}{\left(u\right)}}\right)}{\left(b,v\right)}}\\
			&=\lambdaa{\left(a,u\right)}{\left(\alphaa{-1}{\betaa{-1}{\alphaa{}{v}{\lambdaa{\bar{b}}{\lambdaam{\bar{b}}{\alphaa{-1}{v}{\left(a\right)}}}}}{\betaa{}{b}{\lambdaa{\bar{v}}{\lambdaam{\bar{b}}{\betaa{-1}{b}{\left(u\right)}}}}}}{\rhoo{\alphaa{}{\bar{v}}{\lambdaam{\bar{b}}{\alphaa{-1}{v}{\left(a\right)}}}}{\left(b\right)}},\betaa{-1}{\bar{a}}{\rhoo{\betaa{}{\bar{b}}{\lambdaam{\bar{v}}{\betaa{-1}{b}{\left(u\right)}}}}{\left(v\right)}}\right)}\\
			&=\lambdaa{\left(a,u\right)}{\left(\alphaa{-1}{\bar{u}}{\rhoo{\lambdaam{b}{\left(a\right)}}{\left(b\right)}}, \betaa{-1}{\bar{a}}{\rhoo{\lambdaam{v}{\left(u\right)}}{\left(v\right)}}\right)}\qquad\qquad\qquad\qquad\qquad\mbox{by \eqref{eq:terzo} and \eqref{eq:quarto}}\\
			&=\left(\alphaa{}{u}{\lambdaa{\bar{a}}{\alphaa{-1}{\bar{u}}{\rhoo{\lambdaam{b}{\left(a\right)}}{\left(b\right)}}}},\betaa{}{a}{\lambdaa{\bar{u}}{\betaa{-1}{\bar{a}}{\rhoo{\lambdaam{v}{\left(u\right)}}{\left(v\right)}}}}\right)\\
			&=\left(\lambdaa{a}{\rhoo{\lambdaam{b}{\left(a\right)}}{\left(b\right)}}, \lambdaa{u}{\rhoo{\lambdaam{v}{\left(u\right)}}{\left(v\right)}}\right)\quad\qquad\qquad\qquad\qquad\qquad\qquad\mbox{by \eqref{eq:terzo} and \eqref{eq:quarto}}\\
			&=\left(a\triangleright_{r_S}b, u\triangleright_{r_T}v \right).
		\end{align*}
		This completes the proof.
	\end{proof}
\end{theor}

In other words, \cref{th:shleft} implies \cite[Proposition $1$]{CCoSt19x}. Indeed, that proposition states that the derived solution of the matched product of two left non-degenerate solutions $r_S$ and $r_T$ is exactly the direct product of the derived solutions of $r_S$ and of $r_T$, hence it does not depend on $\alpha$ and $\beta$. 
Recall that the \emph{derived solution} of a left non-degenerate solution r is defined by
\begin{align*}
	r'\left(x,y\right):=\left(y,\lambdaa{y}{\rhoo{\lambdaam{x}{\left(y\right)}}{\left(x\right)}}\right),
\end{align*}
and clearly $r'\left(x,y\right) = \left(y, y\triangleright_r x\right)$, for all $x,y \in X$.

\section*{Funding}
This work was partially supported by the Dipartimento di Matematica e Fisica ``Ennio De Giorgi", Università del Salento. 
The authors are members of GNSAGA (INdAM).

\bibliographystyle{abbrv}
\bibliography{bibliography}

\end{document}